\documentclass[
graybox]{svmult}

\usepackage{
amsmath,
amssymb}

\usepackage{mathptmx}       
\usepackage{helvet}         
\usepackage{courier}        
\usepackage{type1cm}        
\usepackage{makeidx}         
\usepackage{graphicx}        
\usepackage{multicol}        
\usepackage[bottom]{footmisc}

\input cyracc.def
 
\makeindex             
\newcommand{\Real}{\mathbb R}
\newcommand{\RPlus}{\Real^{+}}
\newcommand{\plus}{+}

\newcommand{\Nat}{\mathbb N}

\newcommand{\cF}{\mathcal{F}}
\newcommand{\cT}{\mathcal{T}}
\newcommand{\bP}{\mathbf{P}}
\newcommand{\bE}{\mathbf{E}}
\newcommand{\bC}{\mathbf{C}}

\newcommand{\Mes}{\text{\it Mes}}
\newcommand{\ti}{\widetilde}
\newcommand{\Ind}{{\mbox{\sc I}}}
\newcommand{\FF}{\mbox{F}}
\newcommand{\ff}{\mbox{f}}
\newcommand{\HH}{\mbox{H}}
\newcommand{\bbN}{\mathbb N}
\newcommand{\bbI}{\mathbb I}
\newcommand{\bbB}{\mathbb B}
\newcommand{\fraka}{\mathfrak a}
\newcommand{\frakb}{\mathfrak b}
\newcommand{\frakz}{\mathfrak z}
\newcommand{\frakA}{\mathfrak A}
\newcommand{\frakB}{\mathfrak B}
\newcommand{\frakC}{\mathfrak C}
\newcommand{\Mhat}{\widehat{M}}
\newcommand{\what}{\hat{w}}
\newcommand{\ghat}{\hat{g}}
\newcommand{\tiw}{\tilde{w}}
\newcommand{\tig}{\tilde{g}}

\newcommand{\bary}{\bar{y}}

\newcommand{\card}{\text{card}}

\newcommand{\mtilde}{\widetilde{m}}
\newcommand{\vecM}{\overrightarrow{M}}
\newcommand{\vecN}{\overrightarrow{N}}
\newcommand{\vecx}{\overrightarrow{x}}
\newcommand{\vecm}{\overrightarrow{m}}
\newcommand*{\esssup}{\mathop{\mathrm{ess\,sup}}\displaylimits}

\def\rightmaltese{\protect\vspace*{-1ex}
\begin{flushright}\(\maltese\)\end{flushright}}

\def\rightboxBLACK{\protect\vspace*{-.2ex}
\begin{flushright}\(\blacksquare\)\end{flushright}
}
\def\rightboxWHITE{\protect\vspace*{-.2ex}
\begin{flushright}\(\Box\)\end{flushright}
}
\newenvironment{Proof}{{\sc Proof.\ }}{\rightboxWHITE}
\newenvironment{ProofRemark}{{\sc Proof of remark.\ }}{\rightboxBLACK}
\newenvironment{ProofLemma}{{\sc Proof of lemma.\ }}{\rightmaltese}

\title*{Anglers' fishing problem}
\titlerunning{Anglers' fishing problem}
\author{Anna Karpowicz and Krzysztof Szajowski}
\authorrunning{A.~Karpowicz \and K. Szajowski}
\institute{Anna Karpowicz \at Bank Zachodni WBK, Rynek 9/11, 50-950 Wroc{\l}aw, Poland\\ 
\email{a.m.karpowicz@gmail.com}  \and Krzysztof Szajowski \at Institute of Mathematics and Computer Sci., Wybrze\.ze Wyspia\'nskiego 27, 50-370 Wroc{\l}aw, Poland\\ 
\email{Krzysztof.Szajowski@pwr.wroc.pl}}

\begin{document}

%
%
\maketitle

\abstract*{The considered model will be formulated as related to "the fishing problem" even if the other applications of it are much more obvious. The angler goes fishing. He uses various techniques and he has at most two fishing rods. He buys a fishing ticket for a fixed time. The fishes are caught with the use of different methods according to the renewal processes. The fishes' value and the inter arrival times are given by the sequences of independent, identically distributed (i.i.d.) random variables with the known distribution functions. It forms the marked renewal--reward process. The angler's  measure of satisfaction is given by the difference between the utility function, depending on the value of the fishes caught, and the cost function connected with the time of fishing. In this way, the angler's relative opinion about the methods of fishing is modelled. The angler's aim is to have as much satisfaction as possible and additionally he has to leave the lake before a fixed moment. Therefore his goal is to find two optimal stopping times in order to maximize his satisfaction. At the first  moment, he changes the technique of fishing \emph{e.g.} by excluding one rod and intensifying on the rest. Next, he decides when he should stop the expedition. These stopping times have to be shorter than the fixed time of fishing.  
The dynamic programming methods have been used to find these two optimal stopping times and to specify the expected satisfaction of the angler at these times.
}

\abstract{The considered model will be formulated as related to "the fishing problem" even if the other applications of it are much more obvious. The angler goes fishing. He uses various techniques and he has at most two fishing rods. He buys a fishing ticket for a fixed time. The fishes are caught with the use of different methods according to the renewal processes. The fishes' value and the inter arrival times are given by the sequences of independent, identically distributed (i.i.d.) random variables with the known distribution functions. It forms the marked renewal--reward process. The angler's  measure of satisfaction is given by the difference between the utility function, depending on the value of the fishes caught, and the cost function connected with the time of fishing. In this way, the angler's relative opinion about the methods of fishing is modelled. The angler's aim is to have as much satisfaction as possible and additionally he has to leave the lake before a fixed moment. Therefore his goal is to find two optimal stopping times in order to maximize his satisfaction. At the first  moment, he changes the technique of fishing \emph{e.g.} by excluding one rod and intensifying on the rest. Next, he decides when he should stop the expedition. These stopping times have to be shorter than the fixed time of fishing.  
The dynamic programming methods have been used to find these two optimal stopping times and to specify the expected satisfaction of the angler at these times.
}



\section{Introduction}
Before we start the analysis of the double optimal stopping problem (cf. idea of multiple stopping for stochastic sequences in Haggstrom~\cite{hag67:double}, Nikolaev~\cite{nik79:obob}) for the marked renewal process related to the angler behavior, let us present the so called "fishing problem". One of the first authors who considered the basic version of this problem was Starr~\cite{sta74:captime} and further generalizations were done by Starr and Woodroofe~\cite{stawoo74:gone}, Starr, Wardrop and Woodroofe~\cite{stawarwoo76:delobs}, Kramer, Starr~\cite{krasta90:sizedep} et al. The detailed review of the papers related to the "fishing problem" was presented by Ferguson~\cite{fer:pois}. The simple formulation of the fishing problem, where the angler changes the fishing place or technique before leaving the fishing place, has been done by Karpowicz~\cite{kar09:double}. \emph{We extend the problem to a more advanced model by taking into account the various techniques of fishing used the same time (the parallel renewal--reward processes or the multivariate renewal--reward process). It is motivated by the natural, more precise models of the known, real applications of the fishing problem. The typical process of software testing consists of checking subroutines. At the beginning many kinds of bugs are being searched. The consecutive stopping times are moments when the expert stops general testing of modules and starts checking the most important, dangerous type of error. Similarly, in proof reading, it is natural to look for typographic and grammar errors at the same time. Next, we are looking for language mistakes.}

As various works are done by different groups of experts, it is natural that we would compete with each other. If in the first period work is meant for one group and the second period needs other experts, then they can be players of a game between them. In this case the proposed solution is to find the Nash equilibrium where strategies of players are the stopping times.     

The applied techniques of modeling and finding the optimal solution are similar to those used in the formulation and solution of the optimal stopping problem for the risk process. Both models are based on the methodology explicated by Boshuizen and Gouweleeuw~\cite{bosgou93:semi}. The background mathematics for further reading are monographs by Br{\'e}maud~\cite{bre81:poiproc}, Davis~\cite{dav:marmo} and Shiryaev~\cite{shi78:osr}. The optimal stopping problems for the risk process are subject of consideration in papers by Jensen~\cite{jen97:risk}, Ferenstein and Sieroci\'nski~\cite{fersie97:risk}, Muciek~\cite{mucsza07:risk}. A similar problem for the risk process having disruption (\emph{i.e.} when the probability structure of the considered process is changed at one moment $\theta$) has been analyzed by Ferenstein and Pasternak--Winiarski~\cite{ferwin10:disordered}. The model of the last paper brings to mind the change of fishing methods considered here, however it should be made by a decision maker, not the type of the environment.

The following two sections usher details of the model. It is proper to emphasize that the slight modification of the background assumption by adopting multivariate tools (two rods) and the possible control of their numbers in use extort a different structure of the base model (the underlining process, sets of strategies -- admissible filtrations and stopping times). This modified structure allows the introduction of a new kind of knowledge selection which consequently leads to a game model of the anglers' expedition problem in the section~\ref{competfishing} and~\ref{anglersGAME}. After a quite general formulation a version of the problem for a detailed solution will be chosen. However, the solution is presented as the scalable procedure dependent on parameters which depends on various circumstances. It is not difficult to adopt a solution to wide range of natural cases.      

\subsection{Single Angler's expedition}
The angler goes fishing. He buys a fishing ticket for a fixed time $t_0$ which gives him the right to use at most two rods.  The total cost of fishing depends on real time of each equipment usage and the number of rods used simultaneously. He starts fishing with two rods up to the moment $s$. The effect on each rod can be modelled by the renewal processes $\{N_i(t),t\geq 0\}$, where $N_i(t)$
\index{renewal--reward processes!$N_i(t)$--the number of fishes caught on the rod $i$ to the moment $t$} is the number of fishes caught on the rod $i$, $i \in \frakA:=\{1,2\}$ during the time $t$. Let us combine them together to the marked renewal process. The usage of the $i$-th rod by the time $t$ generates cost $c_i:[0,t_0]\rightarrow \Re$ (when the rod is used simultaneously with other rods it will be denoted by the index dependent on the set of rods, \emph{e.g.}  $\fraka$,  $c^{\fraka}_{i}$) and the reward represented by \emph{i.i.d.} random variables $X^{\{i\}}_{1},X^{\{i\}}_{2},\dots$ 
\index{renewal--reward processes!$X^{\{i\}}_{k}$--the value of the $k$-th fish cached on the $i$-th rod} (\emph{the value of the fishes caught on the $i$-th rod}) with cumulative distribution function $H_i$\footnote{The following convention is used in all the paper: $\vecx=(x_1,x_2,\ldots,x_s)$ for the ordered collection of the elements $\{x_i\}_{i=1}^s$ }. 
The streams of two kinds of fishes are mutually independent and they are independent of the sequence of random moments when the fishes have been caught. The $2$-vector process $\vecN(t)=(N_1(t),N_2(t))$, $t\geq 0$,
\index{renewal--reward processes!$\vecN(t)$--the $2$-dimensional renewal process} can be represented also by a sequence of random variables $T_n$ taking values in $[0,\infty]$ such that
\begin{equation}
\begin{array}{ccc}
&T_0=0,&\\
T_n<\infty&\Rightarrow& T_n<T_{n+1},
\end{array}
\label{pointprocess1}
\end{equation}
\index{renewal--reward processes!$T_n$--$n$-th jump moment}
for $n\in \bbN$, and a sequence of $\frakA$-valued random variables $\frakz_n$
\index{renewal--reward processes!$\frakz_n$--the index of $n$-th jump}
for $n\in \bbN\cup \{0\}$ (see Br\'emaud~\cite{bre81:poiproc} Ch. II, Jacobsen~\cite{jac06:point}). The random variable $T_{n}$ denotes the moment of catching the $n$-th fish ($T_{0}=0$) of any kind and the random variable $\frakz_n$ indicates to which kind the $n$-th fish belongs. The processes $N_i(t)$ can be defined by the sequence $\{(T_n,\frakz_n)\}_{n=0}^\infty$ as:
\begin{equation}
N_i(t)=\sum_{n=1}^\infty \bbI_{\{T_n\leq t\}}\bbI_{\{\frakz_n=i\}}.
\label{point2point}
\end{equation}
Both the $2$-variate process $\vecN(t)$ and the double sequence $\{(T_n,\frakz_n)\}_{n=0}^\infty$ are called \emph{$2$-variate renewal process}. The optimal stopping problems for the compound risk process based on \emph{$2$-variate renewal process} was considered by Szajowski~\cite{sza10:2vector}. 

Let us define, for $i\in \frakA$ and $k\in\bbN$, the sequence 
\begin{equation}
\begin{array}{lcl}
n^{\{i\}}_{0}&=&0,\\
n^{\{i\}}_{k+1}&=&\inf\{n>n^{\{i\}}_{k}:\frakz_n=i\}
\end{array}
\label{index1}
\end{equation}
\index{renewal--reward processes!$n^{\{i\}}_{k}$--the index of $k$-th jump of $i$-th type}%
and put $T^{\{i\}}_{k}=T_{n^{\{i\}}_{k}}$
\index{renewal--reward processes!$T^{\{i\}}_{k}$--$k$-th jump time of the $i$-th type}. 
Let us define random variables $S^{\{i\}}_{n}=T^{\{i\}}_{n}-T^{\{i\}}_{n-1}$
\index{renewal--reward processes!$S^{\{i\}}_{n}$--$n$-th holding time of the $i$-th type}
and assume that they are \emph{i.i.d.} with continuous, cumulative distribution function $F_i(t)=\bP(S^{\{i\}}_{n}\leq t)$
\index{renewal--reward processes!$F_i(t)$--the distribution function of the holding times of the $i$-th type} 
and the conditional distribution function $F_i^s(t)=\bP(S^{\{i\}}_{n}\leq t|S^{\{i\}}_{n}\geq s)$. 
In the section~\ref{PPrep} the alternative representation of the $2$-variate renewal process will be proposed. There is also mild extension of the model in which the stream of events after some moment changes to another stream of events.  
\begin{remark}
In various procedures it is needed to localize the events in a group of the renewal processes. Let $\frakC$ be the set of indices related to such a group. The sequence $\{n^\frakC_k\}_{k=0}^\infty$ such that $n^\frakC_0=0$, $n^\frakC_{k+1}:=\inf\{n>n^\frakC_{k}:\frakz_n\in\frakC \}$ has an obvious meaning. 

Analogously, $n^\frakC(t):=\inf\{n:\! T_n>t,\! \frakz_n\in\frakC\}$. 
\end{remark}

Let $i,j\in\frakA$. The angler's satisfaction measure (\emph{the net reward}) at the period $\fraka$ from the rod $i$ is the difference between the utility function $g^\fraka_{i}: [0,\infty)^2\times\frakA\times \Re^{+} \rightarrow [0,G^\fraka_{i}]$ 
which can be interpreted as the reward from the $i$-th rod when the last success was on rod $j$ and, additionally,  it is dependent on the value of the fishes caught, the moment of results' evaluation, and the cost function $c^\fraka_{i}: [0,t_0] \rightarrow [0,C^\fraka_{i}]$ 
\index{pay-off functions!$c_{i}$, $c^\fraka_{i}$, $c^\fraka$--the cumulative cost of usage the $i$-th rod at the period $\fraka$} reflecting  the cost of duration of the angler's expedition. We assume that $g^\fraka_{i}$ and $c^\fraka_{i}$ are continuous and bounded, additionally $c^\fraka_{i}$ are differentiable. Each fishing method evaluation is based on different utility functions and cost functions. In this way, the angler's relative opinion about them is modelled. 

The angler can change his method of fishing at the moment $s$ and decide to use only one rod. It could be one of the rods used up to the moment $s$ or the other one. Event though the rod used after $s$ is the one chosen from the ones used before $s$ its effectiveness could be different before and after $s$. Following these arguments, the mathematical model of catching fishes, and their value after $s$, could (and in practice should) be different from those for the rods used before $s$. The reason for reduction of the number of rods could be their better effectiveness.  The value of the fishes which have been caught up to time $t$, if the change of the fishing technology took place at the time $s$, is given by
\begin{equation*}
		M_t^s=\sum_{i\in\frakA}\sum_{n=1}^{N_i(s\wedge t)}X^{\{i\}}_{n} + \sum_{n=1}^{N_3((t-s)^+)}X^{\{3\}}_{n}=M_{s\wedge t}+\sum_{n=1}^{N_3((t-s)^{+})}X^{\{3\}}_{n},
\end{equation*}
\index{renewal--reward processes!$M_{t}$ ($M_t^s$)--the renewal-reward process at moment $t$ (with change of a structure at moment $s$)}
where $M^{\{i\}}_{t}=\sum_{n=1}^{N_i(t)}X^{\{i\}}_{n}$,
\index{renewal--reward processes!$M^{\{i\}}_{t}$--the renewal-reward process at moment $t$ related to the rod $i$-th}
and $M_t=\sum_{i=1}^2M^{\{i\}}_{t}$
We denote $\vecM_t=(M^{\{1\}}_{t},M^{\{2\}}_{t})$. Let $Z(s,t)$ denote the angler's pay-off for stopping at time $t$ (the end of the expedition) if the change of the fishing method took place at time $s$. If the effect of extending the expedition after $s$ is described by $g^\frakb_{j}:{\Re^{+}}^2\times\frakA\times [0,t_0]\times \Re\times [0,t_0]\rightarrow [0,G^\frakb_{j}]$, $j\in\frakB$,
\index{pay-off functions!$g^\frakb_{j}(\vecm,i,s,\mtilde,t)$--the reward function after the change of the fishing methods when the state of the renewal-reward processes at $s$ has been $\vecm$ and the final state of the renewal-reward process at $t(\geq s)$ has been $\mtilde$}
minus the additional cost of time $c^\frakb_j(\cdot)$, 
\index{pay-off functions!$c^\frakb_{j}(t)$--the cumulative costs of fishing after the change of fishing method using method $j$}
where $c^\frakb_{j}: [0,t_0] \rightarrow [0,C^\frakb_{j}]$ 
\index{pay-off functions!$C^\frakb_j$--the bounds of the costs} (when $\card(\frakB)=1$ then index $j$ will be abandoned, also $c^\frakb=\sum_{j\in\frakB}c^\frakb_{j}$ will be used, which will be adequate).  
The payoff can be expressed as: 
	\begin{equation}\label{AK Z 1}
		Z(s,t)=\left\{
\begin{array}{ll}
g^\fraka(\vecM_t,\frakz_{N(t)},t)-c^\fraka(t)               &\mbox{ if } t < s \leq t_0, \\
g^\fraka(\vecM_s,\frakz_{N(s)},s)-c^\fraka(s) &\\
+g^\frakb(\vecM_s,\frakz_{N(s)},s,M_t^s,t)-c^\frakb(t-s)&\mbox{ if  $s \leq t \leq t_0$,}\\
 				- C &\mbox{ if }  t_0<t.
\end{array}\right.
\end{equation} 
where the function $c^\fraka(t)$, $g^\fraka(\vecm,i,t)$
\index{pay-off functions!$g^\fraka(\vecm,j,t)$, ($g^\fraka_{i}(\vecm,j,t)$)--the utility of fishes gotten to the moment $t$ (at the $i$-th rod) when the last catch was at the $j$-th rod and the state of the renewal-reward process is $\vecm$} 
and the constant $C$ can be taken as follows: $c^\fraka(t)=\sum_{i=1}^2c^\fraka_{i}(t)$, $g^\fraka(\vecM_s,j,t)=\sum_{i=1}^2g^\fraka_{i}(\vecM_{t},j,t)$, $C=C^\fraka_{1}+C^\fraka_{2}+C^\frakb$. After moment $s$ the modelling process is the renewal--reward one with the stream of \emph{i.i.d.} random variables $X^{\{3\}}_n$ at the moments $T^{\{3\}}_n$ (i.e. appearing according to the renewal process $N_3(t)$). With the notation 
$w^\frakb(\vecm,i,s,\ti m,t)=w^\fraka(\vecm,i,s)+g^\frakb(\vecm,i,s, \ti m,t)-c^\frakb(t-s)$ and $w^\fraka(\vecm,i,t)=g^\fraka(\vecm,i,t)-c^\fraka(t)$, formula~(\ref{AK Z 1}) is reduced to: 
\begin{equation*}
Z(s,t)=Z^{\{\frakz_{N(t)}\}}(s,t)\Ind_{\{t < s \leq t_0\}}+Z^{\{\frakz_{N(s)}\}}(s,t)\Ind_{\{s \leq t \}},
\end{equation*}
where 
	\begin{equation*}
		Z^{\{i\}}(s,t)=\Ind_{\{t < s \leq t_0\}}w^\fraka(\vecM_t,i,t)+\Ind_{\{s \leq t \leq t_0\}}w^\frakb(\vecM_s,i,s,M_t^s,t)-\Ind_{\{t_0<t\}}C.
	\end{equation*}

\subsection{\label{competfishing}The competitive fishing}
When \emph{the methods of fishing are operating by separated anglers} then the stopping random field can be built based on the structure of the marked renewal--reward process as a model of the competitive expedition results. One possible definition of pay-off is based on the assumption that each player has his own account related to the exploration of the fishery. The states of the accounts  depend on who forces the first stop for changing technique, under which circumstances and what techniques they choose. The first stopping moment, the minimum of stopping moments chosen by the players,  is after the moment of the event (catching fish) $T_n$ by the rod $\frakz_n$ and the reward functions depend on the type of fishing which gives recent fish (i.e. $j$, where $j=\frakz_n$). The player's pay-off  $w^\fraka_{i}(\vecm,j,t)=g^\fraka_{i}(\vecm,j,t)-c^\fraka_{i}(t)$
\index{pay-off functions!$w^\fraka_{i}(\vecm,j,t)$--the $i$-th player's pay-off at moment $t$ when the stop has been made by the $j$-th and the state of the renewal-reward process $\vecm$ }. 
The part of the pay-off which depends on the second chosen moment, which stops the expedition, is different for the player who forces the change of fishing methods (the leader) by himself and the other the opponent. The leader is the responsible angler for determining the expedition deadline.

Lets assume for a while that the $i$-th player, $i=1,2$, will take the rod of the opponent and gives his rod to him. \emph{It is not a crucial assumption anyway and the method of fishing after the change can be different from both available before the considered moment.} The method of treatment of the case without this assumption will be explained later (see page~\pageref{SecondFishingProgram}), when the behavior of the player in the second part of the expedition will be formulated.  Define the function
\[
\tiw^\frakb_{i}(\vecm,j,s,k,\mtilde,t)=\tiw^\fraka_{i}(\vecm,j,s)+\tig^\frakb_{i}(\vecm,j,s,k,\mtilde,t)-c^\frakb(t-s)
\]
\index{pay-off functions!$\tiw^\fraka_{i}(\vecm,j,s,k,\mtilde,t)$--the pay-off of the angler $i$-th at moment $t$, when his change of fishing method to $k\in\frakB$ has been forced by the angler $j$ at $s(\leq t)$ and the state of the renewal-reward process $\vecm$ }
for $j\in\frakA$, $k\in\frakB$, where $j$ is the rod by which the fish had been caught just before the moment of the first stop and $k$ is the technique used by $i$-th player after the change (the denotation $-k$ is used for a complimentary rod or player who has decided, which is appropriate). It describes the case when the player deciding to change the method chooses the perspective technique of fishing as the first one. Presumably he will explore the best methods with improvements and the second angler will use the rod which is not used by the leader. The pay-off of the players, when $i$-th is the one who forces the first stop, has the following form:
\begin{eqnarray}
\label{Zij(s,t)1}   Z_{i}(j,s,t)&=&\Ind_{\{t \leq s \leq t_0\}}\tig^\fraka_{i}(\vecM_{t},j,t)+\Ind_{\{s < t \leq t_0\}}\tiw^\frakb_{i}(\vecM_s,i,s,-i,M_t^s,t)-\Ind_{\{t_0<t\}}C\\
\label{Zij(s,t)2}		Z_{-i}(j,s,t)&=&\Ind_{\{t \leq s \leq t_0\}}\tig^\fraka_{-i}(\vecM_{t},j,t)+\Ind_{\{s < t \leq t_0\}}\tiw^\frakb_{-i}(\vecM_s,i,s,i,M_t^s,t)-\Ind_{\{t_0<t\}}C.
\end{eqnarray}
\index{pay-off functions!$Z_{i}(j,s,t)$--the pay--off process of the anglers, when the first stop has been forced by $i$-th one}%
In the above pay--offs it is assumed that the final stop can be declared at any moment. The change of techniques declaration each player makes just after an event at his rod (the catching fish at his rod) as long as on the opponent's rod there is no event. The details of the strategy sets and the solution concept are formulated in the further parts of the paper.


The extension considered here is motivated by the natural, more precise models of the known real applications of the fishing problem. The typical process of software testing consists of checking subroutines. Various types of bugs can be discovered. Each problem with subroutines generates the cost of a bug removal and increases the value of the software. It depends on the types of the bug found. The preliminary testing requires various types of experts. The stable version of subroutines can be kept by less educated computer scientists. The consecutive stopping times are moments when the expert of the defined class stops testing one module and the another tester starts checking. Similarly as in the proof reading. 

\section{The optimization problem and a two person game\label{formulation_game}}
\subsection{\label{PPrep}Filtrations and Markov moments}
Let the sequences of pairs $\{(T_n,\frakz_n)\}_{n=0}^\infty$ be 2-variate renewal process ($\frakA$-marked renewal process) defined on $(\Omega,\cF,\bP)$. According to the denotation of the previous section there are three renewal processes $\{T^{\{i\}}_{n}\}_{n=0}^\infty$, $i=1,2,3$, and 
denoted by $T_n=T^{\{\frakz_n\}}_{N_{\frakz_n}(T_n+)}$.  There are also three renewal--rewarded processes $\{(T^{\{i\}}_{n},X^{\{i\}}_{n})\}_{n=0}^\infty$, $i=1,2,3$
\index{renewal--reward processes!$(T^{\{i\}}_{n},X^{\{i\}}_{n})$--the renewal--rewarded processes}. 
By convention let us denote $X_n=X^{\{\frakz_n\}}_{N_{\frakz_n}(T_n)}$. The following $\sigma$-field generated by history of the $\frakA$-marked renewal processes are defined
\begin{equation}
 \cF_t=\cF_{t}^{\frakA}=\sigma(X_{0},T_{0},\frakz_{0}\ldots,X_{N(t)},T_{N(t)},\frakz_{N(t)}),
\label{sfild}
\end{equation}
\index{filtrations!$\cF_t$--the filtration generated by the $\frakA$-marked renewal--rewarded process to the moment $t$}
\index{filtrations!$\cF_t$,$\cF_{t}^{\frakA}$--the filtration generated by the $\frakA$-marked renewal--rewarded process to the moment $t$}
for $t\geq 0$. This $\sigma$-field can be defined as 
$$\cF_{t}^{\frakA}=\sigma\{(\vecN(s),X_{N(s)},\frakz_{N(s)}), 0\leq s\leq t,\mbox{ $i\in\frakA$}\}.$$ 
\begin{definition}\label{stoppingtimeset}
Let $\cT$ be a set of stopping times with respect to $\sigma$-fields $\{\cF_{t}\}$, $t\geq 0$, defined by \eqref{sfild}. The restricted sets of stopping times are 
\begin{equation}\label{reststopset1}
\cT_{n,K}=\{\tau\in\cT: \mbox{ $\tau\geq 0$, $T_{n}\leq\tau\leq T_{K}$}\} 
\end{equation}
\index{stopping times!$\cT$, $\cT_{n,K}$--sets of stopping times with respect to $\sigma$-fields $\{\cF_{t}\}$}
for $n\in\bbN$, $n<K$ are subsets of $\cT$. The elements of $\cT_{n,K}$
\index{stopping times!$\cT_{n,K}$--the subset of stopping times $\tau\in \cT$ with respect to the filtration $\{\cF_{t}\}$ such that $T_{n}\leq\tau\leq T_{K}$}
are denoted $\tau_{n,K}$
\index{stopping times!$\tau_{n,K}$--the element of the set $\cT_{n,K}$}.
\end{definition}
The stopping times $\tau\in \cT$ have nice representation which will be helpful in the solution of the optimal stopping problems for the renewal processes (see Br{\'e}maud~\cite{bre81:poiproc}). The crucial role in our subsequent considerations plays such a representation. The following lemma is for the unrestricted stopping times.  
\begin{lemma}\label{AK representation1}
	If $\tau \in \cT$ then there exist \mbox{$R_{n}\in \Mes(\cF_n)$} such that the condition $\tau \wedge T_{n+1}=(T_{n}+R_{n})\wedge T_{n+1}$ on $\{\tau \geq T_{n}\}$ a.s. is fulfilled. 
	\end{lemma} 
Various restrictions in the class of admissible stopping times will change this representation.  Some examples of subclasses of $\cT$ are formulated here  (see Lemma~\ref{AK representation1}). Only a few of them are used in optimization problems investigated in the paper (see page~\pageref{payoff}, Corollary~\ref{corollary_AK representation}). 
   
Let $\cF_{s,t}= \sigma(\cF_{s}^{\frakA},X^{\{3\}}_{0},T^{\{3\}}_{0}, \dots,X^{\{3\}}_{N_3((t-s))^\plus},T^{\{3\}}_{N_3((t-s)^\plus)})$ be the $\sigma$-field generated by all events up to time $t$ if the switch at time $s$ from $2$-variate renewal process to another renewal process took place. For simplicity of notation we set\footnote{For the optimization problem there are two epochs: before the first stop, where there are some pay-offs, the model of stream of events, and after the first stop, when there are other pay-offs and different streams of events. In section~\ref{AK second} this  will be emphasized, by adopting adequate denotations.} $\cF^{\{i\}}_{n}:=\cF_{T^{\{i\}}_{n}}$
\index{filtrations!$\cF^{\{i\}}_{n}$--the short denotation of $\cF_{T^{\{i\}}_{n}}$}, 
$\cF_{n}:=\cF_{T_{n}}$, $\cF^s_{n}:=\cF_{s,T^{\{3\}}_{n}}$.  
Let $\Mes(\cF_n)$ ($\Mes(\cF^{\{i\}}_{n})$) denote the set of non-negative and $\cF_n$ ($\cF^{\{i\}}_{n}$)-measurable random variables. From now on, $\cT$ and $\cT^s$ 
stands for the sets of stopping times with respect to $\sigma$-fields $\cF_s$ and $\{\cF_{s,t}, 0\leq s \leq t\}$, respectively. Furthermore, we can define for $n\in \Nat$ and $n\leq K$ the sets
\begin{enumerate}
\item $\cT^{\{i\}}_{n,K}=\{\tau \in \cT:\tau \geq 0,\ T^{\{i\}}_{n}\leq \tau \leq T_{K}\}$;\index{stopping times!$\cT^{\{i\}}_{n,K}$--the stopping times bounded by $T^{\{i\}}_{n}$ and $T_{K}$}
\item $\cT^{\{i\}}_{n}=\{\tau \in \cT:\tau\geq T^{\{i\}}_{n}\}$;\index{stopping times!$\cT^{\{i\}}_{n}$--the stopping times bounded by $T^{\{i\}}_{n}$}
\item $\bar{\cT}^{\{i,\frakA^{\{-i\}}\}}_{n,K}=\{\tau \in \cT:\tau \geq 0,\ T^{\{i\}}_{n}\leq \tau \leq T_{K},\!\; \forall_{k}\tau\notin [T^{\frakA^{-i}}_k,T^{\frakA^{-i}}_{k+1}\vee T^{\{i\}}_{n^{\{i\}}(T^{\frakA^{-i}}_k)}]\}$ where $\frakA^{\{-i\}}:=\frakA\setminus\{i\}$, $T^{\frakA^{-i}}_k:=\min_{\{j\in\frakA^{\{-i\}}\}}\{T^{\{j\}}_{n^{\{j\}}(T^{\{i\}}_k)}\}$;
\item $\bar{\cT}^{\{i\}}_{n}=\{\tau \in \cT:\tau\geq T^{\{i\}}_{n}\!\; , \forall_{k}\tau\notin [T^{\frakA^{-i}}_k,T^{\frakA^{-i}}_{k+1}\vee T^{\{i\}}_{n^{\{i\}}(T^{\frakA^{-i}}_k)}]\}$;
\item $\cT_{n,K}^s=\{\tau \in \cT^s: 0 \leq s \leq \tau,\ T^{\{3\}}_{n}\leq \tau \leq T_{K}\}$.
\end{enumerate}
The stopping times $\tau\in \cT^{\{i\}}$ and $\tau\in \bar{\cT}^{\{i\}}$ can also be represented in the way shown in Lemma~\ref{AK representation1}. 
\begin{lemma}\label{AK representation23} Let the index $i\in\frakA$ be chosen and fixed.
\begin{enumerate}
\item  For every $\tau \in \cT^{\{i\}}$ and $n\in\Nat$ there exist \mbox{$R^{\{i\}}_{n}\in \Mes(\cF^{\{i\}}_n)$} such that $\tau \wedge T^{\{i\}}_{n+1}=(T^{\{i\}}_{n}+R^{\{i\}}_{n})\wedge T^{\{i\}}_{n+1}$ on $\{\tau^{\{i\}} \geq T^{\{i\}}_{n}\}$ a.s. is fulfilled. 		
\item 	If $\tau \in \bar{\cT}^{\{i\}}$ and $n\in\Nat$ there exist \mbox{$R^{\{i\}}_{n}\in \Mes(\cF^{\{i\}}_n)$} such that the condition $\tau \wedge T^{\{i\}}_{n+1}=(T^{\{i\}}_{n}+R^{\{i\}}_{n})\wedge T^{\{i\}}_{n+1}$ on $\{\tau \geq T^{\{i\}}_{n}\}$ a.s.
is fulfilled. 		
\end{enumerate}
	\end{lemma} 

Obviously the angler wants to have as much satisfaction as possible and he has to leave the lake before the fixed moment. Therefore, his goal is to find two optimal stopping times $\tau^{\fraka^*}$\index{stopping times!$\tau^{\fraka^*}$--the optimal moment of the first decision} and $\tau^{\frakb^*}$\index{stopping times!$\tau^{\frakb^*}$--the optimal moment of the second decision} so that the expected gain is maximized 
	\begin{equation}
		\bE Z(\tau^{\fraka^*},\tau^{\frakb^*})=\sup_{\tau^\fraka \in \cT}\sup_{\tau^\frakb \in \cT^{\tau^\fraka}}\bE Z(\tau^\fraka,\tau^\frakb),
	\end{equation}
where $\tau^{\fraka^*}$ corresponds to the moment, when he eventually should change the two rods to the more effective one and $\tau^{\frakb^*}$, when he should stop fishing. These stopping moments should appear before the fixed time of fishing $t_0$. The process $Z(s,t)$ is piecewise-deterministic and belongs to the class of semi-Markov processes. The optimal stopping of similar semi-Markov processes was studied by Boshuizen and Gouweleeuw~\cite{bosgou93:semi} and the multivariate point process by Boshuizen~\cite{bos94:multi}. Here the structure of multivariate processes is discovered and their importance for the model is shown.   
We use the dynamic programming methods to find these two optimal stopping times and to specify the expected satisfaction of the angler. The way of the solution is similar to the methods used by Karpowicz and Szajowski~\cite{karsza07:risk}, Karpowicz~\cite{kar09:double} and Szajowski~\cite{sza10:2vector}. Let us first observe that by the properties of conditional expectation we have
	\begin{eqnarray*}\label{AK double stop b}
		\bE Z(\tau^{\fraka^*},\tau^{\frakb^*})
	 	=\sup_{\tau^\fraka \in \cT}\bE\{\bE\left[Z(\tau^\fraka,\tau^{\frakb^*})|\cF_{\tau^\fraka}\right]\} \nonumber 
	 	= \sup_{\tau^\fraka \in \cT}\bE J(\tau^\fraka),
	\end{eqnarray*}
where
	\begin{equation}\label{AK one stop}
		J(s)=\bE\left[Z(s,\tau^{\frakb^*})|\cF_s\right]=\esssup_{\tau^\frakb \in \cT^{s}}\bE\left[Z(s,\tau^\frakb)|\cF_s\right].
	\end{equation}
Therefore, in order to find $\tau^{\fraka^*}$ and $\tau^{\frakb^*}$, we have to calculate $J(s)$ first. The process $J(s)$ corresponds to the value of the revenue function in one stopping problem if the observation starts at the moment $s$.

\subsection{\label{anglersGAME}Anglers' games}
Based on the consideration of the section~\ref{competfishing} a version of competitive fishing is formulated here. There are two anglers, each using one method of fishing at the beginning of an expedition and an additional fishing period after a certain moment by another method up to the moment chosen by a certain rule. The random field which is the model of payoffs in such a case is given by \eqref{Zij(s,t)1} and \eqref{Zij(s,t)2}. The final segment starts at the moment when one of the anglers wants it. Let $\tau_i\in\bar\cT^{\{i\}}$, $i=\frakA$, be the strategies of the players to stop individual fishing period and switch to the time segment  which is stopped at moment $\sigma$ determined by one angler (let us call them \emph{a leader}). The payoffs of the players are
\begin{eqnarray}\label{payoff}
\psi_i(\tau_1,\tau_2)&=&Z_{i}(\frakz_{N(\tau_1\wedge\tau_2)},\tau_1\wedge\tau_2,\sigma^{\tau_1\wedge\tau_2})\Ind_{\{\tau_1\ne\tau_2\}}\\
\nonumber&&+Z_{i}(\frakz_{N(\tau_1\wedge\tau_2)}\wedge\frakz_{N(\tau_1\wedge\tau_2)},\tau_1\wedge\tau_2,\sigma^{\tau_1\wedge\tau_2})\Ind_{\{\tau_1=\tau_2\}}.
\end{eqnarray}
The aim is to find a pair $(\tau_1^\star,\tau_2^\star)$ of stopping times such that for $i\in \{1,2\}$ we have

\begin{equation}\label{SecondFishingProgram}
\bE\psi_i(\tau_i^\star,\tau_{-i}^\star) \geq \bE\psi_i(\tau_i,\tau_{-i}^\star).
\end{equation}

The optimization problem of the angler and the game between two anglers will involve the construction of the optimal second stopping moment.

\section{Construction of the optimal second stopping time\label{AK second}}
In this section, we will find the solution of one stopping problem defined by~(\ref{AK one stop}). We will first solve the problem for the fixed number of fishes caught, next we will consider the case with the infinite stream  of fishes caught. In this section we fix $s$ - the moment when the change took place and $m=M_s$ - the mass of the fishes at the time $s$. Taking into account various models of fishing after the first stop it is needed to admit various models of stream of events. Assume that the moments of successive fishes catching after the first stop are $T^{\{3\}}_{n}$ and the times between the events are i.i.d. with continuous, cumulative distribution function $\FF(t)$ \index{renewal--reward processes!$\FF(t)$, $\ff(t)$--the distribution and density functions of the holding times after the change of fishing method} 
with the density function $\ff(t)$
and the fishes value represented by \emph{i.i.d.} random variables with distribution function $\HH(t)$\index{renewal--reward processes!$\HH(t)$--the distribution function of the rewards after the change of fishing method} 
(for conveniences this part of expedition is modelled by the renewal process denoted $(T^{\{3\}}_{n},X^{\{3\}}_{n})$).     

\subsection{Fixed number of fishes caught}\label{second fixed}
In this subsection we are looking for optimal stopping time  ${\tau^{\frakb^*}_{0,K}}
:={\tau^\frakb_{K}}^*$\index{stopping times!${\tau^\frakb_{0,K}}^*$, ${\tau^\frakb_{K}}^*$--the second optimal stopping time in a restricted problem}
	\begin{equation}
  	\bE\left[Z(s,{\tau^\frakb_{K}}^*)|\cF_{s}\right]=\esssup_{\tau^\frakb_{K} \in \cT^{s}_{0,K}}\bE\left[Z(s,\tau^\frakb_{K})|\cF_{s}\right],
	\end{equation}
where $s\geq 0$ is a fixed time when the position was changed and $K$ is the maximum number of fishes which can be caught. Let us define
	\begin{equation}\label{Gamma_{n,K}^s}
  	       \Gamma_{n,K}^s=\esssup_{\tau^\frakb_{n,K} \in \cT_{n,K}^s}\bE\left[Z(s,\tau^\frakb_{n,K})|\cF^{s}_{n}\right]%
  	         	=\bE\left[Z(s,{\tau^{\frakb^*}_{n,K}})|\cF^{s}_{n}\right],\ \ \ n=K,\dots,1,0
	\end{equation}
and observe that $\Gamma_{K,K}^s=Z(s,T^{\{3\}}_{K})$. In the subsequent considerations we will use the representation of stopping time formulated in Lemma~\ref{AK representation1} and \ref{AK representation23}. The exact form of the stopping strategies are given in the following corollary.
	
\begin{corollary}\label{corollary_AK representation}
Let $i\in\frakA$.	If $\tau^\fraka \in \cT^{\{i\}}$, $\tau^\frakb \in \cT^s$, then there exist \mbox{$R^\fraka_{n}\in \Mes(\cF^{\{i\}}_n)$} and \mbox{$R^\frakb_{n}\in \Mes(\cF^{s}_{n})$} respectively, such that for conditions $\tau^\fraka \wedge T^{\{i\}}_{n+1}=(T^{\{i\}}_{n}+R^\fraka_{n})\wedge T^{\{i\}}_{n+1}$ on $\{\tau^\fraka \geq T^{\{i\}}_{n}\}$ a.s. and $\tau^\frakb \wedge T^{\{3\}}_{n+1}=(T^{\{3\}}_{n}+R^\fraka_{n})\wedge T^{\{3\}}_{n+1}$ on $\{\tau^\fraka \geq s\wedge T^{\{3\}}_{n}\}$ a.s. are valid.  
	\end{corollary}  

Now we can derive the dynamic programming equations satisfied by $\Gamma_{n,K}^s$. To simplify the notation we can write $M_t=M_t^s$ for $t\leq s$, $\Mhat^{\{1\}}_n=M_{T^{1}_{n}}$, $M_n^s=M_{T^{\{3\}}_{n}}^s$ and $\bar{\FF}_i=1-\FF_i$. The payoff functions are simplified here to $\ghat^\fraka(m)=g^\fraka(m_1,m_2,i,t)\Ind_{\{m_1+m_2=m\}}(m_1,m_2)$, $\ghat^\frakb(m)=g^\frakb(m_1,m_2,i,s,\ti m,t)\Ind_{\{\ti m-m_1-m_2=m\}}$

	\begin{lemma}\label{Gamma recursion_lemma}
		Let $s \geq 0$ be the moment of changing fishery. For $n= K-1, K-2, \dots, 0$
		\begin{equation}\label{optrec}	
			\begin{array}{rcl}
 				\Gamma_{K,K}^s&=& Z(s,T^{\{3\}}_{K}),\\  
    		\Gamma_{n,K}^s&=& \esssup_{R^\frakb_{n}\in \Mes(\cF^{s}_{n})} \vartheta_{n,K}(M_s,s,M_{n}^s,T^{\{3\}}_{n},R^\frakb_{n}) \mbox{ a.s.,}    		    \end{array}
     \end{equation}		
    	 where	
    	\begin{eqnarray*}	
 			\vartheta_{n,K}(m,s,\ti m,t,r) &=& \Ind_{\{t \leq t_0 \}}\bigg\{\bar{\FF}(r)
 								[\Ind_{\{r \leq t_0-t\}}\what^\frakb(m,s,\ti m,t+r)-C\Ind_{\{r>t_0-t\}}]\\
    														&+& E\left[\Ind_{\{ S^{\{3\}}_{n+1}\leq r \}}\Gamma_{n+1,K}^s|\cF^{s}_{n}\right]\bigg\}
    		-C\Ind_{\{t > t_0 \}}    
  		\end{eqnarray*} 
and there exists ${R^\frakb_{n}}^\star\in \Mes(\cF^{s}_{n})$ such that
\begin{equation}\label{optRex}
\Gamma_{n,K}^s=\vartheta_{n,K}(M_s,s,M_{n}^s,T^{\{3\}}_{n},{R^\frakb_{n}}^\star)\mbox{ a.s.,}
\end{equation} 
\begin{equation}\label{tau_2star}
		\tau^{\frakb^*}_{n,K}=\left\{
			\begin{array}{ll}
				{\tau^{\frakb^*}_{n+1,K}} & \mbox{if ${R^\frakb_{n}}^* \geq S^{\{3\}}_{n+1}$,}\\
				T^{\{3\}}_{n}+{R^\frakb_{n}}^* & \mbox{if ${R^\frakb_{n}}^* < S^{\{3\}}_{n+1}$,}\\
		   \end{array}\right.
	\end{equation} 
${\tau^{\frakb^*}_{K,K}}=T^{\{3\}}_{K}$ and $\what^\frakb(m,s,\ti m,t)=\what^\fraka(m,s)+\ghat^\frakb(\ti m - m)-c^\frakb(t-s)$ where $\what^\fraka(m,t)=\ghat^\fraka(m)-c^\fraka(t)$. 		
\end{lemma}

\begin{remark}\label{OptimaSecondTime}
   Let $\{{R^{\frakb^*}_{n}}\}_{n=1}^K$, ${R^{\frakb^*}_{K}}=0$, be a sequence of $\cF^{s}_{n}$--measurable random variables, $n=1,2,\ldots,K$, and ${\eta_{n,K}^{\star s}}=K\wedge \inf\{i\geq n: {R^\frakb_i}^\star<S^{\{3\}}_{i+1}\}$. Then  
  		$\Gamma_{n,K}^s = \bE\left[Z(s,{\tau^{\frakb^*}_{n,K}})|\cF^{s}_{n}\right]$ for $n \leq K-1$, where ${\tau^{\frakb^*}_{n,K}}=T_{\eta_{n,K}^{\star s}}+R^{\frakb^\star}_{{\eta_{n,K}^{\star s}}}$.
\end{remark}

\begin{ProofRemark} \ref{OptimaSecondTime}.
It is a consequence of an optimal choice ${R^\frakb}^\star_n$ in (\ref{optrec}).
\end{ProofRemark}

\begin{ProofLemma}\ref{Gamma recursion_lemma}
First observe that the form of the $\Gamma_{n,K}^s$ for the case $T^{\{3\}}_{n} > t_0$ is obvious from (\ref{AK Z 1}) and (\ref{Gamma_{n,K}^s}).  Let us assume (\ref{optrec}) and (\ref{optRex}) for $n+1,n+2,\ldots,K$. For any $\tau \in \cT_{n,K}^s$ (i.e. $\tau\geq T^{\{3\}}_{n}$ we have 
		$\{\tau<T^{\{3\}}_{n+1}\}=\{\tau\wedge T^{\{3\}}_{n+1} < T^{\{3\}}_{n+1}\} =\{T^{\{3\}}_{n}+R^\frakb_{n}< T^{\{3\}}_{n+1}\}. $
It implies 
	\begin{equation}\label{taus2inc}
		\{\tau<T^{\{3\}}_{n+1}\}=\{S^{\{3\}}_{n+1}>R^\frakb_{n}\}, \ \ \ 
		\{\tau \geq T^{\{3\}}_{n+1}\}=\{S^{\{3\}}_{n+1} \leq R^\frakb_{n}\}.
	\end{equation}
Suppose that $T^{\{3\}}_{K-1}\leq t_0$ and take any $\tau^\frakb_{K-1,K} \in \cT_{K-1,K}^s$.  
According to (\ref{taus2inc}) and the properties of conditional expectation
	\begin{eqnarray*}
		\bE\left[Z(s,\tau)|\cF^{s}_{n}\right]
		&=& \bE\left[\Ind_{\{S^{\{3\}}_{n+1} \leq R^\frakb_{n} \}}\bE[Z(s,\tau \vee T^{\{3\}}_{n+1})|\cF^{s}_{n+1}]|\cF_n\right]\\
    &+& \bE\left[\Ind_{\{S^{\{3\}}_{n+1} > R^\frakb_{n}\}}Z(s,\tau \wedge T^{\{3\}}_{n+1})|\cF^{s}_{n}\right]\\
    &=& \Ind_{\{R^\frakb_{n}\leq t_0-T_n \}}\bar{\FF}(R_n)\what^\frakb(M_s,s,M_{n}^s,T^{\{3\}}_{n}+R^\frakb_{n})\\
    &+& \bE\left[\Ind_{\{S^{\{3\}}_{n+1} \leq R^\frakb_{n}\}}\bE[Z(s,\tau \vee T^{\{3\}}_{n+1})|\cF^{s}_{n+1}|\cF^{s}_{n}\right].   
	\end{eqnarray*} 
Let $\sigma\in\cT^\frakb_{n+1}$. For every $\tau\in\cT_n$ we have
\begin{equation*}
		\tau=\left\{
			\begin{array}{ll}
				\sigma & \mbox{if $R^\frakb_{n} \geq S^{\{3\}}_{n+1}$,}\\
				T^{\{3\}}_{n}+R^\frakb_{n}  & \mbox{if $R^\frakb_{n} < S^{\{3\}}_{n+1}$.}\\
		   \end{array}\right.
	\end{equation*} 
We have 
\begin{eqnarray*}
\bE[Z(s,\tau)|\cF^{s}_{n}]&=&\bE\left[\Ind_{\{S^{\{3\}}_{n+1} \leq R^\frakb_{n} \}}\bE[Z(s,\sigma)|\cF^{s}_{n+1}]|\cF_n\right]\\
    &+& \Ind_{\{R^\frakb_{n}\leq t_0-T_n \}}\bar{\FF}(R^\frakb_{n})\what^\frakb(M_s,s,M_{n}^s,T^{\{3\}}_{n}+R^\frakb_{n})\\
    &\leq& \sup_{R\in \Mes(\cF^{s}_{n})}\{\bE\left[\Ind_{\{S^{\{3\}}_{n+1} \leq R \}}\Gamma^s_{n+1,K}|\cF_n\right]\\
    &+& \Ind_{\{R\leq t_0-T_n \}}\bar{\FF}(R)\what^\frakb(M_s,s,M_{n}^s,T^{\{3\}}_{n}+R)\}=\bE[Z(s,\tau^\star_{n,K}|\cF^s_{n}]
\end{eqnarray*}
It follows $\sup_{\tau\in \cT^{s}_{n}}\bE[Z(s,\tau)|\cF^s_{n}]\leq \bE[Z(s,\tau^\star_{n,K}|\cF^s_{n}]\leq \sup_{\tau\in \cT^\frakb_{n}}\bE[Z(s,\tau)|\cF^s_{n}]$ where the last inequality is because $\tau^\star_{n,K}\in \cT^{s}_{n,K}$.
We apply the induction hypothesis, which completes the proof.
\end{ProofLemma} 
\begin{lemma}\label{AK Gamma 2} 
		$\Gamma_{n,K}^s=\gamma_{K-n}^{s,M_s}(M_{n}^s,T^{\{3\}}_{n})$ for $n=K,\dots,0$, where the sequence of functions 							$\gamma_j^{s,m}$ is given recursively as follows:
			\begin{eqnarray}\label{gamma recursive a}
    		\gamma_{0}^{s,m}(\ti m,t)&=&\Ind_{\{t\leq t_0\}}\what^\frakb(m,s,\ti m,t)-C\Ind_{\{t > t_0\}},\nonumber \\
    		\gamma_{j}^{s,m}(\ti m,t)&=&\Ind_{\{t\leq t_0\}}\sup_{ r\geq 0}\kappa^\frakb_{\gamma_{j-1}^{s,m}}(m,s,\ti m,t,r)
    																	-  C\Ind_{\{t > t_0\}},
    	\end{eqnarray}
    where		
    	\begin{eqnarray*}																
    		\kappa^\frakb_{\delta}(m,s,\ti m,t,r)&=&
    		\bar{\FF}(r)[\Ind_{\{r \leq t_0-t \}}\what^\frakb(m,s,\ti m,t+r)-C\Ind_{\{r > t_0-t\}}]\nonumber\\
    		&+& \int_0^r d\FF(z) \int_0^\infty \delta(\ti m+x,t+z)d\HH(x).
  		\end{eqnarray*}
	\end{lemma}

\begin{ProofLemma}\ref{AK Gamma 2}.
Since the case for $t > t_0$ is obvious let us assume that $T^{\{3\}}_{n}\leq t_0$ for $n \in \{0,\dots,K-1\}$. Let us notice that according to Lemma~\ref{Gamma recursion_lemma} we obtain   
		$\Gamma_{K,K}^s 
									 = \gamma_{0}^{s,M_s}(M_{K}^s,T^{\{3\}}_{K})$,
thus the proposition is satisfied for $n=K$. Let $n=K-1$ then Lemma~\ref{Gamma recursion_lemma} and the induction hypothesis leads to
	\begin{eqnarray*}
    \Gamma_{K-1,K}^s&=& \esssup_{R^\frakb_{K-1}\in \Mes(\cF_{s,K-1})}\bigg\{
    \bar{\FF}(R^\frakb_{K-1})[\Ind_{\{R^\frakb_{K-1} \leq t_0-T^{\{3\}}_{K-1} \}}\what^\frakb(M_s,s,M_{K-1}^s,T^{\{3\}}_{K-1}+R^\frakb_{K-1})\\
    &-&C\Ind_{\{R^\frakb_{K-1}>t_0-T^{\{3\}}_{K-1}\}}]
    +\bE\left[\Ind_{\{ S^{\{3\}}_{K}\leq R^\frakb_{K-1} \}}\gamma_{0}^{s,M_s}(M_{K}^s,T^{\{3\}}_{K})|\cF_{s,K-1}\right]\bigg\} \ a.s.,
	\end{eqnarray*}
where $M_{K}^s=M_{K-1}^s+X^{\{3\}}_{K}$, $T^{\{3\}}_{K}=T^{\{3\}}_{K-1}+S^{\{3\}}_{K}$ and the random variables $X^{\{3\}}_{K}$ and $S^{\{3\}}_{K}$ are independent of $\cF_{s,K-1}$. Moreover $R^\frakb_{K-1}$, $M^s_{K-1}$ and 
$T^{\{3\}}_{K-1}$ are $\cF_{s,K-1}$-measurable. It follows 
	\begin{eqnarray*}
    \Gamma_{K-1,K}^s
    & = & \esssup_{R^\frakb_{K-1}\in \Mes(\cF_{s,K-1})}\bigg\{\bar{\FF}(R^\frakb_{K-1})
   					[\Ind_{\{R^\frakb_{K-1} \leq t_0-T^{\{3\}}_{K-1} \}}\what^\frakb(M_s,s,M_{K-1}^s,T^{\{3\}}_{K-1}+R^\frakb_{K-1})\\
    &-&C\Ind_{\{R^\frakb_{K-1}>t_0-T^{\{3\}}_{K-1}\}}]
    			+ \int_0^{R^\frakb_{K-1}} d\FF(z) \int_0^\infty \gamma_{0}^{s,M_s}(M_{K-1}^s+x,T^{\{3\}}_{K-1}+z)d\HH(x)\bigg\}\\
    & = & \gamma^{s,M_s}_1(M^s_{K-1},T^{\{3\}}_{K-1}) \ a.s.
	\end{eqnarray*}
Let $n \in \{1,\dots,K-1\}$ and suppose that $\Gamma_{n,K}^s=\gamma_{K-n}^{s,M_s}(M_{n}^s,T^{\{3\}}_{n})$. Similarly like before,  we conclude by Lemma~\ref{Gamma recursion_lemma} and induction hypothesis that
  \begin{eqnarray*}
    \Gamma_{n-1,K}^s&=& \esssup_{R^\frakb_{n-1}\in \Mes(\cF^{s}_{n-1})}\bigg\{
    \bar{\FF}(R^\frakb_{n-1})[\Ind_{\{R^\frakb_{n-1} \leq t_0-T^{\{3\}}_{n-1} \}}\what^\frakb(M_s,s,M_{n-1}^s,T^{\{3\}}_{n-1}+R^\frakb_{n-1})\\
    &-&C\Ind_{\{R^\frakb_{n-1}>t_0-T^{\{3\}}_{n-1}\}}]
    +\int_0^{R^\frakb_{n-1}}d\FF(s)\int_0^\infty \gamma_{K-n}^{s,M_s}(M_{n-1}^s+x,T^{\{3\}}_{n-1}+s)d\HH(x)\bigg\} \ a.s.
	\end{eqnarray*}
therefore $\Gamma_{n-1,K}^s=\gamma_{K-(n-1)}^{s,M_s}(M_{n-1}^s,T^{\{3\}}_{n-1})$.
\end{ProofLemma}

From now on we will use $\alpha_i$ to denote the hazard rate of the distribution $F_i$ (i.e. $\alpha_i=f_i/\bar{\FF}_i$) and to shorten notation we set $\Delta^\cdot(a)=\bE\left[\ghat^\cdot(a+X^{\{i\}})-\ghat^\cdot(a)\right]$, where $\cdot$ can be $\fraka$ or $\frakb$.
\begin{remark}\label{gamma recursive 2}
		The sequence of functions $\gamma_j^{s,m}$ can be expressed as:
			\begin{eqnarray*}
    		\gamma_{0}^{s,m}(\ti m,t)&=&\Ind_{\{t\leq t_0\}}\what^\frakb(m,s,\ti m,t)-C\Ind_{\{t > t_0\}},\\
    		\gamma_{j}^{s,m}(\ti m,t)&=&\Ind_{\{t\leq t_0\}}\bigg\{ \what^\frakb(m,s,\ti m,t)+y^\frakb_{j}(\ti m-m,t-s,t_0-t) \bigg\}
    	-C\Ind_{\{t > t_0\}} \nonumber
		\end{eqnarray*} 
		and  $y^\frakb_{j}(a,b,c)$ is given recursively as follows 
		\begin{eqnarray*}
    	y^\frakb_{0}(a,b,c)&=& 0  \\
    	y^\frakb_{j}(a,b,c)&=&\max_{ 0\leq r \leq c}\phi^\frakb_{y^\frakb_{j-1}}(a,b,c,r),
 		\end{eqnarray*}  
 		where
 		$ \phi^\frakb_{\delta}(a,b,c,r)=  \int_0^r \bar{\FF}(z)
    													\{\alpha_2(z)\left[\Delta^\frakb(a)+\bE\delta(a+X^{\{3\}},b+z,c-z)\right]
    							- {c^\frakb}^\prime(b+z)\}dz. $
\end{remark}
	
\begin{ProofRemark}\ref{gamma recursive 2}
Clearly
	\begin{equation*}
		\int_0^r d\FF(s) \int_0^\infty \gamma_{j-1}^{s,m}(\ti m+x,t+s)d\HH(x)
			=\bE\left[\Ind_{\{S^{\{3\}} \leq r \}}\gamma_{j-1}^{s,m}(\ti m+X^{\{3\}},t+S^{\{3\}})\right],
	\end{equation*}
where $S^{\{3\}}$ has \emph{c.d.f.} $\FF$ and $X^{\{3\}}$ has c.d.f. $\HH$. Since $\FF$ is continuous and $\kappa^\frakb_{\gamma_{j-1}^{s,m}}(m,s,\ti m,t,r)$ is bounded and continuous for $t\in \RPlus \setminus \{t_0\}$, the supremum in (\ref{gamma recursive a}) can be changed into maximum. Let $r>t_0-t$ then
	\begin{eqnarray*}
    \kappa^\frakb_{\gamma_{j-1}^{s,m}}(m,s,\ti m,t,r)
    											&=& \bE\left[\Ind_{\{S^{\{3\}} \leq t_0-t \}}
    											\gamma_{j-1}^{s,m}(\ti m+X^{\{3\}},t+S^{\{3\}})\right]-C\bar{\FF}(t_0-t)\\
    											&\leq & \bE\left[\Ind_{\{S^{\{3\}} \leq t_0-t \}}
    											\gamma_{j-1}^{s,m}(\ti m+X^{\{3\}},t+S^{\{3\}})\right]+\bar{\FF}(t_0-t)\what^\frakb(m,s,\ti m,t_0)\\
    											&=& \kappa^\frakb_{\gamma_{j-1}^{s,m}}(m,s,\ti m,t,t_0-t).
	\end{eqnarray*}
The above calculations cause that 
    		$\gamma_{j}^{s,m}(\ti m,t)=\Ind_{\{t\leq t_0\}}\max_{ 0\leq r \leq t_0-t}\varphi_{j}(m,s,\ti m,t,r)
    														 - C\Ind_{\{t > t_0\}}$, 
where 
		$\varphi_{j}(m,s,\ti m,t,r)=\bar{\FF}(r)\what^\frakb(m,s,\ti m,t+r)
		+ \bE\left[\Ind_{\{S^{\{3\}} \leq r \}}\gamma_{j-1}^{s,m}(\ti m+X^{\{3\}},t+S^{\{3\}})\right]$.
Obviously for $S^{\{3\}} \leq r$ and $r\leq t_0-t$ we have $S^{\{3\}} \leq t_0$ therefore we can consider the cases $t \leq t_0$ and $t > t_0$ separately. Let $t \leq t_0$ then $\gamma_0^{s,m}(\ti m,t)=\what^\frakb(m,s,\ti m,t)$ and the hypothesis is true for $j=0$. The task is now to calculate $\gamma_{j+1}^{s,m}(\ti m,t)$ given $\gamma_{j}^{s,m}(\cdot,\cdot)$. The induction hypothesis implies that for $t \leq t_0$ 
	\begin{eqnarray*}
	\varphi_{j+1}(m,s,\ti m,t,r)&=& \bar{\FF}(r)\what^\frakb(m,s,\ti m,t+r)
    												 						+ \bE\left[\Ind_{\{S^{\{3\}} \leq r \}}\gamma_{j}^{s,m}(\ti m+X^{\{3\}},t+S^{\{3\}})\right]\\
    &=& \ghat^\fraka(m)-c^\fraka(s)+\bar{\FF}(r)\left[\ghat^\frakb(\ti m-m)-c^\frakb(t-s+r)\right]\\
    												 &+& \int_0^r
    																		\ff(z)\{\bE \ghat^\frakb(\ti m-m+X^{\{3\}})-c^\frakb(t-s+z)\\
    												 &+& \bE y^\frakb_{j}(\ti m-m+X^{\{3\}},t-s+z,t_0-t-z)\}dz. 		
	\end{eqnarray*}
It is clear that for any $a$ and $b$ 
	\begin{eqnarray*}
		\bar{\FF}(r)\left[\ghat^\frakb(a)-c^\frakb(b+r)\right]
						&=&\ghat^\frakb(a)-c^\frakb(b)\\
						&&\mbox{}-\int_0^r\{\ff(z)\left[\ghat^\frakb(a)-c^\frakb(b+z)\right]+\bar{\FF}(z){c^\frakb}^\prime(b+z)\}dz,  
	\end{eqnarray*} 
therefore 
	\begin{eqnarray*}
		\varphi_{j+1}(m,s,\ti m,t,r) &=& \what^\frakb(m,s,\ti m, t)
 				+ \int_0^r \bar{\FF}(z)\{ \alpha_2(z)[\Delta^\frakb(\ti m-m)\\ 
 			 &+& \bE y^\frakb_{j}(\ti m-m+X^{\{3\}},t-s+z,t_0-t-z)] - {c^\frakb}^\prime(t-s+z)\}dz, 									
   \end{eqnarray*}
which proves the theorem.  
The case for $t>t_0$ is trivial.
\end{ProofRemark}

Following the methods of Ferenstein and Sieroci\'nski~\cite{fersie97:risk}, we find the second optimal stopping time. Let $\bbB=\frakB([0,\infty)\times[0,t_0]\times[0,t_0])$ be the space of all bounded, continuous functions with the norm \mbox{$\left\|\delta\right\|=\sup_{a,b,c}|\delta(a,b,c)|$}. It is easy to check that $\bbB$ with the norm supremum is complete space. The operator $\Phi^\frakb: \bbB \rightarrow \bbB$ is defined by 
	\begin{equation}\label{AK Phi_2}
		(\Phi^\frakb\delta)(a,b,c)=\max_{0\leq r\leq c}\phi^\frakb_{\delta}(a,b,c,r).
	\end{equation}
Let us observe that $y^\frakb_{j}(a,b,c)=(\Phi^\frakb y^\frakb_{j-1})(a,b,c)$. Remark~\ref{gamma recursive 2} now implies that	there exists a function ${r^\frakb}^*_{j}(a,b,c)$ such that $y^\frakb_{j}(a,b,c)=\phi^\frakb_{y^\frakb_{j-1}}(a,b,c,{r^\frakb}^*_{j}(a,b,c))$ and this gives 
		\begin{eqnarray*}\label{AK r 2}
			\gamma_{j}^{s,m}(\ti m,t)&=&\Ind_{\{t\leq t_0\}}\bigg\{ \what^\frakb(m,s,\ti m,t)\\
			&&\mbox{}+\phi^\frakb_{y^\frakb_{j-1}}(\ti m-m,t-s,t_0-t,{r^{\frakb^*}_{j}}(\ti m-m,t-s,t_0-t))\bigg\}
    	-C\Ind_{\{t > t_0\}}.	
		\end{eqnarray*}
The consequence of the foregoing considerations is the theorem, which determines optimal stopping
times ${\tau^{\frakb^*}_{n,K}}$ in the following manner:
	\begin{theorem}\label{AK solution 2}
	Let ${R^\frakb_{i}}^*  = {r^{\frakb^*}_{K-i}}(M_{i}^s-M_s,T^{\{3\}}_{i}-s,t_0-T^{\{3\}}_{i})$ for $i=0,1,\dots,K$ moreover
			$\eta^s_{n,K} = K\wedge \inf \{i\geq n: {R^\frakb_{i}}^*< S^{\{3\}}_{i+1}\}$, then 
	the stopping time ${\tau^{\frakb^*}_{n,K}}=T^{\{3\}}_{\eta^s_{n,K}}+ R^{\frakb^*}_{\eta^s_{n,K}}$ 
is optimal in the class $\cT^s_{n,K}$ and $\Gamma_{n,K}^s=\bE\left[Z(s,{\tau^{\frakb^*}_{n,K}})|\cF^{s}_{n}\right]$.
\end{theorem}
\subsection{Infinite number of fishes caught}
The task is now to find the function $J(s)$ and stopping time $\tau^{\frakb^*}$, which is optimal in class $\cT^s$. In order to get the solution of one stopping problem for infinite number of fishes caught it is necessary to put the restriction $\FF(t_0)<1$.  
	\begin{lemma}\label{contraction 2}
	If $\FF(t_0)<1$ then the operator $\Phi^\frakb:\bbB \rightarrow \bbB$ defined by~(\ref{AK Phi_2}) is a contraction.
	\end{lemma}
\begin{ProofLemma} \ref{contraction 2}.
Let $\delta_i \in \bbB$ assuming that $i \in \{1,2\}$. There exists $\rho_i $ such that $(\Phi^\frakb\delta_i)(a,b,c)=\phi^\frakb_{\delta_i}(a,b,c,\rho_i)$. We thus get
	\begin{eqnarray*}
		(\Phi^\frakb \delta_1)(a,b,c)-(\Phi^\frakb \delta_2)(a,b,c)
				&=&  \phi^\frakb_{\delta_1}(a,b,c,\rho_1)-\phi^\frakb_{\delta_2}(a,b,c,\rho_2)\\
		&\leq& \phi^\frakb_{\delta_1}(a,b,c,\rho_1)-\phi^\frakb_{\delta_2}(a,b,c,\rho_1)\\
		&=&  \int_0^{\rho_1}d\FF(z)\int_0^{\infty}
										[\delta_1-\delta_2](a+x,b+z,c-s)d\HH(x)\\
		&\leq& \int_0^{\rho_1}d\FF(z)\int_0^{\infty}
										\sup_{a,b,c}|[\delta_1-\delta_2](a,b,c)|d\HH(x)\\  
		&\leq& \FF(c) \left\| \delta_1-\delta_2 \right\| \leq \FF(t_0) \left\| \delta_1-\delta_2 \right\|
				\leq \bC \left\| \delta_1-\delta_2 \right\|,
	\end{eqnarray*}
where $0\leq \bC < 1$. 
Similarly, like as before,
	$(\Phi^\frakb \delta_2)(a,b,c)-(\Phi^\frakb \delta_1)(a,b,c)
				\leq \bC \left\| \delta_2-\delta_1 \right\|$.
Finally we conclude that	
	$\left\|\Phi^\frakb \delta_1-\Phi^\frakb \delta_2\right\|\leq \bC \left\| \delta_1-\delta_2 \right\|$
which completes the proof.
\end{ProofLemma}

Applying Remark~\ref{gamma recursive 2}, Lemma~\ref{contraction 2} and the fixed point theorem we conclude
\begin{remark}\label{converge 2}
		There exists $y^\frakb \in \bbB$ such that $y^\frakb=\Phi^\frakb y^\frakb$  and $\lim_{K\rightarrow \infty }\|y^\frakb_{K}-y^\frakb\|=0$. 
\end{remark}	
According to the above remark, $y^\frakb$ is the uniform limit of $y^\frakb_{K}$, when $K$ tends to infinity, which implies that $y^\frakb$ is measurable and $\gamma^{s,m}=\lim_{K\rightarrow \infty}\gamma^{s,m}_K$ is given by
	\begin{equation}\label{gamma^{s,m}}	
			\gamma^{s,m}(\ti m,t )
			=\Ind_{\{ t \leq t_0\}}\left[\what^\frakb(m,s,\ti m,t)+y^\frakb(\ti m-m,t-s,t_0-t)\right]-C\Ind_{\{t > t_0\}}.
	\end{equation}		
We can now calculate the optimal strategy and the expected gain after changing the place.
	\begin{theorem}\label{solution 2}
	If $\FF(t_0)<1$ and has the density function $\ff$, then
		\begin{itemize}
  		\item[\rm(i)] for $n \in \Nat$ the limit $\tau^{\frakb^\star}_{n}=\lim_{K\rightarrow \infty}{\tau^{\frakb^*}_{n,K}}\ a.s.$ exists and $\tau^{\frakb^\star}_{n} \leq t_0$ is an optimal stopping rule in the set $\cT^s\cap \{\tau \geq T^{\{3\}}_{n}\}$,
   		\item[\rm(ii)] $ \bE\left[Z(s, \tau^{\frakb^\star}_{n})|\cF^{s}_{n}\right]=\gamma^{s,m}(M_{n}^s,T^{\{3\}}_{n})$ a.s.
		\end{itemize}
	\end{theorem}
\begin{Proof} 
(i) Let us first prove the existence of  $\tau^{\frakb^\star}_{n}$. By definition of $\Gamma_{n,K+1}^s$ we have
	\begin{eqnarray*}
		\Gamma_{n,K+1}^s
										&=& \esssup_{\tau \in \cT_{n,K+1}^s}\bE\left[Z(s,\tau)|\cF^{s}_{n}\right]
								=\esssup_{\tau \in \cT_{n,K}^s}\bE\left[Z(s,\tau)|\cF^{s}_{n}\right] 
								\vee \esssup_{\tau \in \cT_{K,K+1}^s}\bE\left[Z(s,\tau)|\cF^{s}_{n}\right]\\
								&=&\bE\left[Z(s,{\tau^{\frakb^*}_{n,K}})|\cF^{s}_{n}\right] \vee \bE\left[Z(s,\sigma^*)|\cF^{s}_{n}\right]
	\end{eqnarray*}
thus we	observe that $\tau^{\frakb^*}_{n,K+1}$ is equal to $\tau^{\frakb^*}_{n,K}$ or $\sigma^*$, 
where $\tau^{\frakb^*}_{n,K} \in \cT_{n,K}^s$ and $\sigma^* \in \cT_{K,K+1}^s$ respectively. 
It follows that $\tau^{\frakb^*}_{n,K+1}\geq \tau^{\frakb^*}_{n,K}$ which implies that the sequence $\tau^{\frakb^*}_{n,K}$ is nondecreasing with respect to $K$. Moreover  ${R^\frakb_{i}}^* \leq t_0 - T^{\{3\}}_{i}$ for all $i \in \{0, \dots, K\}$ thus  
$\tau^{\frakb^*}_{n,K} \leq t_0$ and therefore $\tau^{\frakb^\star}_{n} \leq t_0$ exists.\newline
Let us now look at the process $\xi^s(t)=(t,M_t^s,V(t))$, where $s$ is fixed and $V(t)=t-T^{\{3\}}_{N_3(t)}$. $\xi^s(t)$ is  Markov process with the state space $[s,t_0]\times [m,\infty) \times [0,\infty)$. In a general case the infinitesimal operator for $\xi^{s}$ is given by
	\begin{eqnarray*}
		A p^{s,m}(t,\ti m, v) &=& \frac{\partial}{\partial t} p^{s,m}(t,\ti m, v) 
																+ \frac{\partial}{\partial v} p^{s,m}(t,\ti m, v)\\
													&+& \alpha_2(v)\bigg\{ \int_{\RPlus} p^{s,m}(t,x,0)d\HH(x)- p^{s,m}(t,\ti m, v) \bigg\},			
	\end{eqnarray*}
where $p^{s,m}(t,\ti m, v) : [0,\infty)\times [0,\infty) \times [0,\infty)\rightarrow \Real$ is continuous, bounded, measurable with bounded left-hand derivatives with respect to $t$ and $v$ (see~\cite{bosgou93:semi} and ~\cite{rolschschteu98}). Let us notice that for $t \geq s$ the process $Z(s,t)$ can be expressed  as 
		$Z(s,t)= p^{s,m}(\xi^{s}(t))$,
where	
	\begin{equation*}
		 p^{s,m}(\xi^{s}(t))=\left\{
			\begin{array}{ll}
 				\ghat^\fraka(M_s)-c^\fraka(s)+\ghat^\frakb(M_t^s-M_s)-c^\frakb(t-s) &\mbox{ if } s \leq t \leq t_0,\\
 				- C &\mbox{ if }  t_0<t.
			\end{array}\right.
	\end{equation*}	 
It follows easily that in our case  $A p^{s,m}(t,\ti m, v)=0$ for $t_0<t$ and 
	\begin{equation}\label{A p^{s,m}}
		A p^{s,m}(t,\ti m, v) =\alpha_2(v)[\bE \ghat^\frakb(\ti m+X^{\{3\}}-m)-\ghat^\frakb(\ti m -m)]-{c^\frakb}\prime(t-s) 	
	\end{equation}
for $s \leq t \leq t_0$. The process $p^{s,m}(\xi^s(t))-p^{s,m}(\xi^s(s))-\int_s^t(A p^{s,m})(\xi^s(z))dz$ is a martingale with respect to $\sigma(\xi^s(z), z\leq t)$ which is the same as $\cF_{s,t}$. This can be found in~\cite{dav:marmo}. Since $\tau^{\frakb^*}_{n,K} \leq t_0 $, applying the Dynkin's formula we obtain
	\begin{equation}\label{Dynkin K}
		\bE\left[p^{s,m}(\xi^s(\tau^{\frakb^*}_{n,K}))|\cF^{s}_{n}\right]-p^{s,m}(\xi^s(T^{\{3\}}_{n}))
		=\bE\left[ \int_{T^{\{3\}}_{n}}^{\tau^{\frakb^*}_{n,K}}(Ap^{s,m})(\xi^s(z))dz|\cF^{s}_{n}\right] \ \ a.s.
	\end{equation}
From (\ref{A p^{s,m}}) we conclude that 
	\begin{eqnarray*}\label{Ap^{s,m}dz}
		\int_{T^{\{3\}}_{n}}^{\tau^{\frakb^*}_{n,K}}(Ap^{s,m})(\xi^s(z))dz&=&
										[\bE \ghat^\frakb(M^s_n+X^{\{3\}}-m)-\ghat^\frakb(M^s_n -m)] \int_{T^{\{3\}}_{n}}^{\tau^{\frakb^*}_{n,K}}\alpha_2(z-T^{\{3\}}_{n})dz \nonumber\\
										&-& \int_{T^{\{3\}}_{n}}^{\tau^{\frakb^*}_{n,K}}{c^\frakb}^\prime(z-s)dz.
	\end{eqnarray*}
Moreover let us check that
	\begin{eqnarray*}
& & \left|\int_{T^{\{3\}}_{n}}^{\tau^{\frakb^*}_{n,K}}\alpha_2(z-T^{\{3\}}_{n})dz \right|\leq
						\frac{1}{\bar \FF(t_0)}\int_{T^{\{3\}}_{n}}^{\tau^{\frakb^*}_{n,K}}\ff(z-T^{\{3\}}_{n})dz
			\leq \frac{1}{\bar \FF(t_0)} < \infty, \\	
		& & \left|\int_{T^{\{3\}}_{n}}^{\tau^{\frakb^*}_{n,K}}{c^\frakb}^\prime(z-s)dz \right|
			= \left|c^\frakb({\tau^{\frakb^*}_{n,K}}-s)-c^\frakb(T^{\{3\}}_{n}-s)\right| < \infty, \\		
		& & \left|\bE \ghat^\frakb(M^s_n+X^{\{3\}}-m)-\ghat^\frakb(M^s_n -m)\right| < \infty,
	\end{eqnarray*}
where the two last inequalities result from the fact that the functions  $\ghat^\frakb$ and $c^\frakb$ are bounded. On account of the above observation we can use the dominated convergence theorem and
	\begin{equation}\label{lim Ap^{s,m}dz}
		\lim_{K\rightarrow \infty} \bE\left[ \int_{T^{\{3\}}_{n}}^{\tau^{\frakb^*}_{n,K}}(Ap^{s,m})(\xi^s(z))dz|\cF^{s}_{n}\right]
											= \bE\left[ \int_{T^{\{3\}}_{n}}^{ \tau^{\frakb^\star}_{n}}(Ap^{s,m})(\xi^s(z))dz|\cF^{s}_{n}\right].
	\end{equation}
Since $\tau^{\frakb^\star}_{n}\leq t_0$ applying the Dynkin's formula to the left side of (\ref{lim Ap^{s,m}dz}) we conclude that
	\begin{equation}\label{Dynkin}
		\bE\left[ \int_{T^{\{3\}}_{n}}^{\tau^{\frakb^\star}_{n}}(Ap^{s,m})(\xi^s(z))dz|\cF^{s}_{n}\right]
			= \bE\left[p^{s,m}(\xi^s(\tau^{\frakb^\star}_{n}))|\cF^{s}_{n}\right]-p^{s,m}(\xi^s(T^{\{3\}}_{n}))\ \ a.s.
	\end{equation}
Combining (\ref{Dynkin K}), (\ref{lim Ap^{s,m}dz}) and (\ref{Dynkin}) we can see that 
	\begin{equation} \label{boundary}
		\lim_{K\rightarrow \infty} \bE\left[p^{s,m}(\xi^s(\tau^{\frakb^*}_{n,K}))|\cF^{s}_{n}\right]=\bE\left[p^{s,m}(\xi^s(\tau^{\frakb^*}_{n}))|\cF^{s}_{n}\right],
	\end{equation}
hence that $\lim_{K\rightarrow \infty} \bE\left[Z(s,\tau^{\frakb^*}_{n,K})|\cF^{s}_{n}\right]=\bE\left[Z(s,\tau^{\frakb^*}_{n})|\cF^{s}_{n}\right]$.
We next prove the optimality of ${\tau^\frakb_{n}}^*$ in the class $\cT^s\cap \{\tau^\frakb_{n} \geq T^{\{3\}}_{n}\}$. Let $\tau \in \cT^s\cap \{\tau^\frakb_{n} \geq T^{\{3\}}_{n}\}$ and it is clear that $\tau \wedge T^{\{3\}}_{K} \in \cT^s_{n,K}$. As $\tau^{\frakb^*}_{n,K}$ is optimal in the class $\cT^s_{n,K}$ we have
	\begin{equation}\label{ineq a}
		\lim_{K\rightarrow \infty}\bE\left[p^{s,m}(\xi^s({\tau^{\frakb^*}_{n,K}}))|\cF^{s}_{n}\right]
			\geq \lim_{K\rightarrow \infty}\bE\left[p^{s,m}(\xi^s(\tau \wedge T^{\{3\}}_{K}))|\cF^{s}_{n}\right].
	\end{equation}
From (\ref{boundary}) and (\ref{ineq a}) we conclude that
$\bE\left[p^{s,m}(\xi^s({\tau^{\frakb^*}_{n}}))|\cF^{s}_{n}\right] \geq \bE\left[p^{s,m}(\xi^s(\tau))|\cF^{s}_{n}\right]$
for any stopping time $\tau \in \cT^s\cap \{\tau \geq T^{\{3\}}_{n}\}$, which implies that ${\tau^{\frakb^*}_{n}}$ is optimal in this class. \newline
(ii) Lemma \ref{AK Gamma 2} and (\ref{boundary}) lead to
		$\bE\left[Z(s,{\tau^\frakb_{n}}^*)|\cF^{s}_{n}\right]=\gamma^{s,M_s}(M_{n}^s,T^{\{3\}}_{n})$.
\end{Proof}
The remainder of this section will be devoted to the proof of the left-hand differentiability of the function $\gamma^{s,m}(m,s)$ with respect to $s$. This property is necessary to construct the first optimal stopping time. First, let us briefly denote $\delta(0,0,c)\in \bbB$ by $\bar \delta(c)$.   
	\begin{lemma}\label{operator derivatives}
	Let $\bar \nu(c)=\Phi^\frakb \bar \delta(c)$, $\bar \delta(c) \in \bbB$ and $\left|\bar \delta'_{+}(c)\right|\leq A_1$ for 
	$c \in [0,t_0)$ then $\left|\bar\nu'_{+}(c)\right|\leq A_2$.   		
	\end{lemma}
\begin{ProofLemma} \ref{operator derivatives}.
	First observe that the derivative $\bar \nu'_{+}(c)$ exists because 
	$\bar \nu(c)=\max_{0\leq r \leq c} \bar \phi^\frakb(c,r)$, where  $\bar \phi^\frakb(c,r)$ is differentiable with respect to 
	$c$ and $r$. 
	Fix $h\in(0,t_0-c)$ and define 
	$\bar\delta_1(c) = \bar\delta(c+h)\in \bbB$ and $\bar\delta_2(c) = \bar\delta(c)\in \bbB$. Obviously, 		
  			$\|\Phi^\frakb  \bar\delta_1-\Phi^\frakb  \bar\delta_2 \| \geq |\Phi^\frakb  \bar\delta_1(c)-\Phi^\frakb  \bar\delta_2(c) |
  			= |\Phi^\frakb  \bar\delta(c+h)-\Phi^\frakb  \bar\delta(c) |$
	and on the other side using Taylor's formula for right-hand derivatives we obtain
		\begin{equation*} 
			 \left\|  \bar\delta_1- \bar\delta_2 \right\| 
			 = \sup_{c}\left|\bar\delta(c+h)- \bar\delta(c) \right|																
			 \leq h \sup_{c}\left|\bar \delta'_{+}(c)\right|+\left| o(h) \right|.
		\end{equation*}
		From the above and Remark \ref{contraction 1} it follows that
		\begin{equation*} \label{Together}
			-\bC \bigg \{ \sup_{c}\left|\bar \delta'_{+}(c)\right|+\frac{\left|o(h)\right| }{h} \bigg\} 
												\leq \frac{\bar \nu(c+h)-\bar \nu(c)}{h} 
											\leq \bC \bigg \{ \sup_{c}\left|\bar \delta'_{+}(c)\right|+\frac{\left|o(h)\right| }{h} \bigg\}
		\end{equation*}
and letting $h \rightarrow 0^+$  gives $\left|\bar \nu'_{+}(c)\right|\leq \bC A_1=A_2$. 
\end{ProofLemma}
The significance of Lemma	\ref{operator derivatives} is such that the function $\bary(t_0-s)$ has bounded left-hand derivative with respect to $s$ for $s \in (0,t_0]$. The important consequence of this fact is the following
	\begin{remark}\label{gamma form}
		The function $\gamma^{s,m}$ can be expressed as 
			$\gamma^{s,m}(m,s)=\Ind_{\{s\leq t_0\}}u(m,s)-C\Ind_{\{s > t_0\}}$, 
		where $u(m,s)=\ghat^\fraka(m)-c^\fraka(s)+\ghat^\frakb(0)-c^\frakb(0)+\bary^\frakb(t_0-s)$ is continuous, bounded, measurable with the bounded left-hand derivatives with respect to $s$.
	\end{remark}
At the end of this section, we determine the conditional value function of the second optimal stopping problem. According to (\ref{AK one stop}), Theorem~\ref{solution 2} and Remark~\ref{gamma form} we have  
	\begin{equation}\label{J(s)}
		J(s)=\bE\left[Z(s,\tau^{\frakb^*})|\cF_s\right]=\gamma^{s,M_s}(M_s,s)\text{ a.s.}
	\end{equation}	
	
\section{Construction of the optimal first stopping time}\label{first stopping}
In this section, we formulate the solution of the double stopping problem. 
On the first epoch of the expedition the admissible strategies (stopping times) depend on the formulation of the problem. For the optimization problem the most natural are the stopping times from $\cT$ (see the relevant problem considered in Szajowski~\cite{sza10:2vector}). However, when the bilateral problem is considered the natural class of admissible strategies depends on who uses the strategy. It should be $\cT^{\{i\}}$ for the $i$-th player. Here the optimization problem with restriction to the strategies from the $\cT^{\{1\}}$ at the first epoch is investigated. 

Let us first notice that the function $u(m,s)$ has a similar properties to the function $\what^\frakb(m,s,\ti m, t)$ and the process $J(s)$ has similar structure to the process $Z(s,t)$. By this observation one can follow the calculations of Section~\ref{AK second} to get $J(s)$. Let us define again
  	$\Gamma_{n,K}=\esssup_{\tau^\fraka \in \cT_{n,K}}\bE\left[J(\tau^\fraka)|\cF_{n}\right],\ n=K,\dots,1,0,$
which fulfills the following representation 	
	\begin{lemma}\label{AK Gamma 1}
	$\Gamma_{n,K}=\gamma_{K-n}(\Mhat^{\{1\}}_n,T^{\{1\}}_{n})$ for $n=K,\dots,0,$ where the sequence of functions $\gamma_j$ can be expressed as:
		\begin{eqnarray*}
    		\gamma_{0}(m,s)&=&\Ind_{\{s\leq t_0\}}u(m,s)-C\Ind_{\{s > t_0\}},\\
    		\gamma_{j}(m,s)&=&\Ind_{\{s\leq t_0\}}\bigg\{u(m,s)+y^\fraka_{j}(m,s,t_0-s) \bigg\}
    	-C\Ind_{\{s > t_0\}} \nonumber
		\end{eqnarray*} 
		and  $y^\fraka_{j}(a,b,c)$ is given recursively as follows: 
		\begin{eqnarray*}
    	y^\fraka_{0}(a,b,c)&=& 0  \\
    	y^\fraka_{j}(a,b,c)&=&\max_{ 0\leq r \leq c}\phi^\fraka_{y^\fraka_{j-1}}(a,b,c,r)
 		\end{eqnarray*}  
 		where
\begin{eqnarray*}
\phi^\fraka_{\delta}(a,b,c,r)&=&  \int_0^r \bar{\FF}_1(z)
    													\left\{\alpha_1(z)\left[\Delta^\fraka(a)+\bE\delta(a+x^{\{1\}},b+z,c-z)\right]\right.\\
    													&&\mbox{}-\left.({\bary^\frakb}{}^\prime_{-}(c-z)+{c^\fraka}^\prime(b+z))\right\}dz. 
\end{eqnarray*}
	\end{lemma}
Lemma \ref{AK Gamma 1} corresponds to the combination of Lemma~\ref{AK Gamma 2} and Remark~\ref{gamma recursive 2} from Subsection~\ref{second fixed}.  Let the operator $\Phi^\fraka: \bbB \rightarrow \bbB$ be defined by 
	\begin{equation}\label{AK Phi_1}
		(\Phi^\fraka\delta)(a,b,c)=\max_{0\leq r\leq c}\phi^\fraka_{\delta}(a,b,c,r).
	\end{equation}
Lemma \ref{AK Gamma 1} implies that	there exists a function $r_{1,j}^*(a,b,c)$ such that 
		\begin{equation*}\label{AK r 1}
			\gamma_{j}(m,s)=\Ind_{\{s\leq t_0\}}\bigg\{ u(m,s)+\phi^\fraka_{y^\fraka_{j-1}}(m,s,t_0-s,r_{1,j}^*(m,s,t_0-s))\bigg\}
    	-C\Ind_{\{s > t_0\}}.	
		\end{equation*}
We can now state the analogue of Theorem \ref{AK solution 2}.
\begin{theorem}\label{AK solution 1}
	Let $R^{\fraka^*}_{i}  = r^{\fraka^*}_{K-i}(M_{i},T^{\{1\}}_{i},t_0-T^{\{1\}}_{i})$ and
			$\eta_{n,K} = K\wedge \inf \{i\geq n: {R^\fraka_{i}}^*< S^{\{1\}}_{i+1}\}$, then 
	 ${\tau^{\fraka^*}_{n,K}}=T^{\{1\}}_{\eta_{n,K}}+ {R^{\fraka^*}_{\eta_{n,K}}}$ 
is optimal in the class $\cT_{n,K}$ and $\Gamma_{n,K}=E\left[J(\tau^{\fraka^*}_{n,K})|\cF_{n}\right]$.
\end{theorem}
The following results may be proved in much the same way as in Section~\ref{AK second}. 
	\begin{lemma}\label{contraction 1}
		If $F_1(t_0)<1$ then the operator $\Phi^\fraka:\bbB \rightarrow \bbB$ defined by (\ref{AK Phi_1}) is a contraction.
	\end{lemma}
		\begin{remark}\label{converge 1}
		There exists $y^\fraka \in \bbB$ such that $y^\fraka=\Phi^\fraka y^\fraka$  and $\lim_{K\rightarrow \infty }\|y^\fraka_{K}-y^\fraka\|=0$. 
	\end{remark}	
The above remark implies that $\gamma=\lim_{K\rightarrow \infty}\gamma_K$ is given by
	\begin{equation}\label{gamma}	
			\gamma(m,s)
			=\Ind_{\{ s \leq t_0\}}\left[u(m,s)+ y^\fraka(m,s,t_0-s)\right]-C\Ind_{\{s > t_0\}}.
	\end{equation}
We can now formulate our main results. 
	\begin{theorem}\label{solution 1}
		If $F_1(t_0)<1$ and has the density function $f_1$, then
		\begin{itemize}
  		\item[\rm(i)] for $n \in \Nat$ the limit $\tau^{\fraka^*}_{n}=\lim_{K\rightarrow \infty}\tau^{\fraka^*}_{n,K}\ a.s.$ exists and $\tau^{\fraka^*}_{n} \leq t_0$ is an optimal stopping rule in the set $\cT\cap \{\tau \geq T^{\{1\}}_{n}\}$,
   		\item[\rm(ii)]$ \bE\left[J(\tau^{\fraka^*}_{n})|\cF_{n}\right]=\gamma(M_{n},T^{\{1\}}_{n})\ \ \ \ a.s.$
		\end{itemize}
	\end{theorem}
\begin{Proof}
The proof follows the same method as in Theorem~\ref{solution 2}. The difference lies in the form of the infinitesimal operator. Define the processes $\xi(s)=(s,M_s,V(s))$ where $V(s)=s-T^{\{1\}}_{N_1(s)}$. Like before $\xi(s)$ is the Markov process with the state space $[0,\infty)\times [0,\infty) \times [0,\infty)$. Notice that $J(s)= p(\xi(s))$ and $p(s,m,v) : [0,t_0]\times [0,\infty) \times [0,\infty)\rightarrow \Real$ continuous, bounded, measurable with the bounded left-hand derivatives with respect to $s$ and $v$. It is easily seen that  
		$A p(s,m, v) =\alpha_1(v)[\bE \ghat^\fraka(m+x^{\{1\}})-\ghat^\fraka(m)]-\left[{\bary^{\frakb^\prime}_{-}}(t_0-s)+{c^\fraka}^\prime(s)\right]$ for $s \leq t_0$. 
The rest of the proof remains the same as in the proof of  Theorem~\ref{solution 2}.
\end{Proof}
Summarizing, the solution of a double stopping problem is given by
	\begin{equation*}
		\bE Z(\tau^{\fraka^*},\tau^{\frakb^*})=\bE J(\tau^{\fraka^*})=\gamma(M_{0},T^{\{1\}}_{0})=\gamma(0,0),
	\end{equation*}
where $\tau^{\fraka^*}$ and $\tau^{\frakb^*}$ are defined according to Theorem~\ref{solution 2} and Theorem~\ref{solution 1} respectively. 
\section{Examples}
The form of the solution results in the fact that it is difficult to calculate the solution in an analytic way. In this section we will present examples of the conditions for which the solution can be calculated exactly. 
	\begin{remark}\label{Ap^{s,m} monotone}
		If the process $\zeta_2(t)=Ap^{s,m}(\xi^s(t))$ has decreasing paths, then the second optimal stopping time is given by
			${\tau^{\frakb^*}_{n}}=\inf \{ t \in \left[T^{\{3\}}_{n},t_0\right] : Ap^{s,m}(\xi^s(t)) \leq 0 \}$
on the other side if $\zeta_2(t)$ has non-decreasing paths, then the second optimal stopping time is equal to $t_0$. 
		\newline
Similarly, if the process $\zeta_1(s)=Ap(\xi(s))$ has decreasing paths, then the first optimal stopping time is given by
			$\tau^{\fraka^*}_{n}=\inf \{ s\in \left[T^{\{1\}}_{n},t_0\right] : Ap(\xi(s)) \leq 0 \}$
on the other side if $\zeta_1(s)$ has non-decreasing paths, then the first optimal stopping time is equal to $t_0$.
	\end{remark}
\begin{Proof} 
From (\ref{Dynkin}) we obtain 
		$\bE\left[Z(s,{\tau^{\frakb^*}_{n}})|\cF^{s}_{n}\right]
			=Z(s,T^{\{3\}}_{n})+\bE\left[ \int_{T^{\{3\}}_{n}}^{{\tau^{\frakb^*}_{n}}}(Ap^{s,m})(\xi^s(z))dz\right]$ a.s.
and the application results of Jensen and Hsu~\cite{jen:hsu} completes the proof. 
\end{Proof}
	\begin{corollary}
	If $S^{\{3\}}$ has exponential distribution with constant hazard rate $\alpha_2$, the function $\ghat^\frakb$ is increasing and
	concave, the cost function $c^\frakb$ is convex and $t_{2,n}=T^{\{3\}}_{n}$, $m_n^s=M_n^s$ then
	\begin{equation}\label{tau^*_{2,n}}
		{\tau^{\frakb^*}_{n}}=\inf \{ t \in \left[t_{2,n},t_0 \right] : \alpha_2[\bE \ghat^\frakb(m_n^s + x^{\{3\}}-m)-\ghat^\frakb(m_n^s -m)]\leq {c^\frakb}^\prime(t-s) \},
	\end{equation}
	where $s$ is the moment of changing the place. Moreover, if $S^{\{1\}}$ has exponential distribution with constant hazard rate $\alpha_1$, $\ghat^\fraka$ is increasing and concave, $c^\fraka$ is convex and $t_{1,n}=T^{\{1\}}_{n}$, $m_n=\Mhat^{\{1\}}_n$ then
	$$\tau^{\fraka^*}_{n}=\inf \{ s \in \left[t_{1,n},t_0 \right] : \alpha_1\left[\bE \ghat^\fraka(m_n+ x^{\{1\}})-\ghat^\fraka(m_n)\right]\leq  {c^\fraka}^\prime(s) \}$$
	\end{corollary}
	\begin{Proof}
	The form of ${\tau^\fraka}^*_{n}$ and ${\tau^\frakb_{n}}^*$ is a consequence of Remark~\ref{Ap^{s,m} monotone}. Let us observe that	by our assumptions $\zeta_2(t)=\alpha_2\Delta^\frakb(M_t^s-m)-{c^\frakb}^\prime(t-s)$ has decreasing paths for $t \in [T^{\{3\}}_{n},T^{\{3\}}_{n+1})$. It suffices to prove that $\zeta_2(T^{\{3\}}_{n})-\zeta_2(T^{{2}-}_{n})=\alpha_2[\Delta^\frakb(M_{n}^s-m)-\Delta^\frakb(M_{n-1}^s-m)]<0$ for all $n\in \Nat$. \\
	It remains to check that ${\bary^{\frakb^\prime}}_{-}(t_0-s)=0$. We can see that ${\tau^\frakb}^*={\tau^\frakb}^*(s)$ is deterministic, which is clear from (\ref{tau^*_{2,n}}). Let us notice that if $s \leq t_0$ then combining (\ref{Dynkin}), (\ref{boundary}) and (\ref{J(s)}) gives
		$\gamma^{s,m}(m,s)=\bE\left[Z(s,\tau^{\frakb^*})|\cF_s\right]=Z(s,s)
																							+\bE\left[ \int_{s}^{ {\tau^\frakb}^*}(Ap^{s,m})(\xi^s(z))dz|\cF_s\right].$			By Remark \ref{gamma form} it follows that 
$${\bary}^\frakb(t_0-s)=\bE\left[ \int_{s}^{ {\tau^\frakb}^*(s)}(Ap^{s,m})(\xi^s(z))dz\right]
		=\int_{s}^{ {\tau^\frakb}^*(s)}\left[\alpha_2\Delta^\frakb(0)-c'_2(z-s)\right]dz
$$ 
and this yields 
		\begin{eqnarray}\label{y'}
			\bary^{\frakb^\prime}_{-}(t_0-s)&=&\int_{s}^{ {\tau^\frakb}^*(s)}c''_2(z-s)dz
			+{{\tau^\frakb}^*}'(s)\left[\alpha_2\Delta^\frakb(0)-c'_2({\tau^\frakb}^*_{2}(s)-s)\right] \\
\nonumber			&&-\left[\alpha_2\Delta^\frakb(0)-c'_2(0)\right]\\
\nonumber			&=& c'_2({\tau^\frakb}^*(s)-s)-c'_2(0)-\left[\alpha_2\Delta^\frakb(0)-c'_2(0)\right]=0.
\end{eqnarray}					
		The rest of proof runs as before.
	\end{Proof}
	\begin{corollary}
	If for $i=1$ and $i=2$ the functions $g_i$ are increasing and convex, $c_i$ are concave and $S^{\{i\}}$ have the exponential distribution with constant hazard rate $\alpha_i$ then $\tau^{\fraka^*}_{n}={\tau^{\frakb^*}_{n}}=t_0$ for $n \in \Nat$.
	\end{corollary}	
	\begin{Proof}
	It is also the straightforward consequence of Remark~\ref{Ap^{s,m} monotone}. It suffices to check that $\bary^{\frakb^\prime}_{-}(t_0-s)$ is non-increasing with respect to $s$. First observe that $\tau^{\frakb^*}(s)=t_0$. Considering (\ref{y'}) it is obvious that $\bary^{\frakb^\prime}_{-}(t_0-s)=\alpha_2\Delta^\frakb(0)-c'_2(t_0-s)$ and this completes the proof.
	\end{Proof}
\section{Conclusions}	
This article presents the solution of the double stopping problem in the "fishing model" for the finite horizon. The analytical properties of the reward function in one stopping problem played the crucial rule in our considerations and allowed us to get the solution for the extended problem of a double stopping.  Let us notice that by repeating considerations from Section~\ref{first stopping} it is easy to generalize our model and the solution to the multiple stopping problem but the notation can be inconvenient. The construction of the equilibrium in the two person non-zero sum problem formulated in the section~\ref{formulation_game} can be reduced to the two double optimal stopping problems in the case when the payoff structure is given by \eqref{Zij(s,t)1}, \eqref{Zij(s,t)2} and \eqref{payoff}. The key assumptions were related to the properties of the distribution functions. Assuming general distributions and the infinite horizon one can get the extensions of the above model.      



\printindex
\end{document}